\documentclass[production,12pt]{jsg}
\usepackage{latexsym,amsmath,amsfonts,amscd,amssymb,enumerate,url}

\newtheorem{thm}{Theorem}

\newtheorem{defi}[thm]{Definition}

\newtheorem*{Rmks}{Remarks}
\newtheorem*{Rmk}{Remark}

\newtheorem*{Xpl}{Example}

\newtheorem*{Pb}{Problem}

\newcommand\al{\alpha}

\newcommand\om{\omega}
\newcommand\sub{\subseteq}

\newcommand{\Z}{\mathbb{Z}}

\newcommand{\R}{\mathbb{R}}

\newcommand\x{\times}

\newcommand\A{\mathcal A}

\newcommand\wt[1]{{\widetilde{#1}}}

\newcommand\UU{\mathcal{U}}

\newcommand{\hhat}[1]{\widehat{#1}}

\newcommand\nn{{\nonumber}}

\newcommand\id{{\operatorname{id}}}

\renewcommand\phi{\varphi}

\newcommand\wrt{w.r.t.~}
\newcommand\Fix{{\operatorname{Fix}}}
\newcommand\Fixc{{\operatorname{Fix}_{\operatorname{c}}}}

\newcommand\CF{{\operatorname{CF}}}
\newcommand\HF{\operatorname{HF}}
\newcommand\dd{\partial}

\newcommand\D{\mathbb{D}}
\newcommand\reff[1]{(\ref{#1})}
\newcommand\Jreg{\mathcal{J}_{\operatorname{reg}}}

\title[Coisotropic Floer homology]{Note on coisotropic Floer homology and leafwise fixed points}

\author{Fabian Ziltener}

\begin{document}

\begin{abstract}
For an adiscal or monotone regular coisotropic submanifold $N$ of a symplectic manifold I define its Floer homology to be the Floer homology of a certain Lagrangian embedding of $N$. Given a Hamiltonian isotopy $\phi=(\phi^t)$ and a suitable almost complex structure, the corresponding Floer chain complex is generated by the $(N,\phi)$-contractible leafwise fixed points. I also outline the construction of a local Floer homology for an arbitrary closed coisotropic submanifold.

Results by Floer and Albers about Lagrangian Floer homology imply lower bounds on the number of leafwise fixed points. This reproduces earlier results of mine.

The first construction also gives rise to a Floer homology for a Boothby-Wang fibration, by applying it to the circle bundle inside the associated complex line bundle. This can be used to show that translated points exist.
\end{abstract}

\maketitle

\section*{Introduction}\label{sec:intro}
Consider a symplectic manifold $(M,\om)$, a coisotropic submanifold $N\sub M$, and a Hamiltonian diffeomorphism $\psi:M\to M$. The isotropic (or characteristic) distribution $TN^\om$ on $N$ gives rise to the isotropic foliation on $N$. A \emph{leafwise fixed point for $\psi$} is a point $x\in N$ for which $\psi(x)$ lies in the leaf through $x$ of this foliation. We denote by $\Fix(\psi,N)$ the set of such points. A fundamental problem in symplectic geometry is the following:
\begin{Pb} Find conditions under which $\Fix(\psi,N)$ is non-empty and find lower bounds on its cardinality.
\end{Pb}
This generalizes the problems of showing that a given Hamiltonian diffeomorphism has a fixed point and that a given Lagrangian submanifold intersects its image under a Hamiltonian diffeomorphism. References for solutions to the general problem are provided in \cite{ZiLeafwise,ZiC0}. 
\begin{Xpl}[translated points] As explained in \cite[p.~97]{SaMorse}, translated points of the time-1-map of a contact isotopy starting at the identity are leafwise fixed points of the Hamiltonian lift of this map to the symplectization.
\end{Xpl}
We denote
\[N_\om:=\big\{\textrm{isotropic leaves of }N\big\}.\]
We call $N$ \emph{regular} (or ``fibering'') iff there exists a smooth manifold structure on $N_\om$, such that the canonical projection from $N$ to $N_\om$ is a smooth submersion.\footnote{Such a structure is unique if it exists. In this case the symplectic quotient of $N$ is well-defined.} %
Let $X$ be a manifold and $h\in C\big([0,1]\x X,M\big)$. We call $h$ a \emph{semistrong $(N,\om)$-homotopy} iff for every $t\in[0,1]$ there is a leaf $F$ of $N$ that contains the image of $\{t\}\x\dd X$ under $h$. 

We denote by $\D$ the closed unit disk in $\R^2$. We call $(N,\om)$ (or simply $N$) \emph{adiscal} iff every map $u\in C\big(\D,M\big)$ that sends $\dd\D=S^1$ to a leaf of $N$, is semistrongly $(N,\om)$-homotopic to a constant map.

The main result of \cite{ZiLeafwise} (Theorem 1.1) implies the following. We denote by $b_i(N)$ the $i$-th $\Z_2$-Betti-number of $N$.
\begin{thm}[leafwise fixed points for adiscal coisotropic]\label{thm:adiscal} Assume that $(M,\om)$ is geometrically bounded, $N$ is closed\footnote{This means compact and without boundary.}, %
regular, and adiscal, and that $(N,\psi)$ is nondegenerate in the sense of \cite[p.~105]{ZiLeafwise}. Then the following estimate holds:
\begin{equation}\label{eq:Fix sum}\big|\Fix(\psi,N)\big|\geq\sum_{i=0}^{\dim N}b_i(N).
\end{equation}
\end{thm}
This bound is sharp if there exists a $\Z_2$-perfect Morse function on $N$, see \cite[Theorem 1.2]{ZiLeafwise}. The idea of the proof of Theorem \ref{thm:adiscal} given in \cite{ZiLeafwise}, is to find a suitable Lagrangian embedding $\iota_N$ of $N$ into a geometrically bounded symplectic manifold, see \reff{eq:iota N N hhat M} below. We then apply Y.~Chekanov's Main Theorem in \cite{Ch}, which implies the result in the Lagrangian case.\footnote{\cite[Theorem 1.1]{ZiLeafwise} is formulated in a more general setting than Theorem \ref{thm:adiscal}. Chekanov's result is needed to deal with that setting, whereas in the setting of Theorem \ref{thm:adiscal} Floer's original article \cite{FlLag} suffices.} %end footnote
The Lagrangian intersection points of the image of $\iota_N$ correspond to leafwise fixed points of $\psi$.

Similarly to Theorem \ref{thm:adiscal}, in \cite{ZiMaslov} for a regular $N$, we defined monotonicity and the minimal Maslov number $m(N)$, and we proved the following result (\cite[Theorem 3]{ZiMaslov}): 

\newpage

\begin{thm}[leafwise fixed points for monotone coisotropic]\label{thm:monotone} Assume that $(M,\om)$ is geometrically bounded or convex at infinity\footnote{\cite[Theorem 1.1]{ZiLeafwise} is stated for the geometrically bounded case, but the proof goes through in the convex at infinity case.}, $N\sub M$ is closed, monotone, and regular, and $(N,\psi)$ is non-degenerate. Then the following estimate holds:
\begin{equation}\label{eq:Fix sum m N}|\Fix(\psi,N)|\geq\sum_{i=\dim N-m(N)+2}^{m(N)-2}b_i(N).\end{equation}
\end{thm}
The idea of the proof of this theorem given in \cite{ZiMaslov}, is to use the same Lagrangian embedding as in the proof of Theorem \ref{thm:adiscal}. We then apply P.~Albers' Main Theorem in \cite{AlLag}, which states Theorem \ref{thm:monotone} in the Lagrangian case.

%on the other hand
Finally, the main result of \cite{ZiC0} (Theorem 1) implies that leafwise fixed points exist for an \emph{arbitrary} closed coisotropic submanifold if the Hamiltonian flow is suitably $C^0$-small. More precisely, it implies the following:
\begin{thm}[leafwise fixed points for $C^0$-close coisotropic]\label{thm:C0} Let $(M,\om)$ be a symplectic manifold and $N\sub M$ be a closed coisotropic submanifold. Then there exists a $C^0$-neighbourhood $\UU\sub C(N,M)$ of the inclusion $N\to M$ with the following property. Let $(\phi^t)_{t\in[0,1]}$ be a Hamiltonian flow on $M$ satisfying $\phi^t|_N\in\UU$, for every $t\in[0,1]$. If $(N,\phi^1)$ is nondegenerate then the following estimate holds:
\begin{equation}\label{eq:Fix sum }\big|\Fix(\phi^1,N)\big|\geq\sum_{i=0}^{\dim N}b_i(N).
\end{equation}
\end{thm}
This result is optimal in the sense that the $C^0$-condition cannot be replaced by Hofer smallness of $\phi^1$. The idea of the proof of Theorem \ref{thm:C0} given in \cite{ZiC0}, is to find a suitable embedding of $N$ as a Lagrangian submanifold $\wt N$ of some symplectic manifold $\wt M$ that is a Weinstein neighbourhood of $\wt N$, see \reff{eq:wt N} below. We then use Weinstein's neighbourhood theorem and the existence of Lagrangian intersection points for the zero-section in the cotangent bundle. These points correspond to leafwise fixed points for $\phi^1$. Since $N$ is not assumed to be regular in Theorem \ref{thm:C0}, in contrast with Theorems \ref{thm:adiscal} and \ref{thm:monotone}, we can construct $\wt M$ only locally around $\wt N$.\\

The point of this note is to reinterpret the proofs of Theorems \ref{thm:adiscal} and \ref{thm:monotone} in terms of a version of Floer homology for an adiscal or monotone regular coisotropic submanifold. I also outline a definition of a local version of Floer homology for an arbitrary closed coisotropic submanifold and use it to reinterpret the proof of Theorem \ref{thm:C0}. Details of the construction of this homology will be carried out elsewhere. For the extreme cases $N=M$ and $N$ Lagrangian, local versions of Floer homology were developed in \cite{FlFixed,OhFloer,OhLocal,CFHW,Po,GGLocal}; see also the book \cite[Chapter 17.2]{OhSympl}.
\footnote{In \cite{AlFloer} a Lagrangian Floer homology was constructed that is ``local'' in a different sense.}

Potentially a (more) global version of coisotropic Floer homology may be defined under a suitable condition on $N$ that is weaker than regularity, so that the $C^0$-condition on $(\phi^t)$ in Theorem \ref{thm:C0} can be relaxed.%
\footnote{This can only work under suitable conditions on $N$. The reason is that by \cite[Theorem 1.1]{GGFrag}, there exists a closed hypersurface $N$ in $\R^{2n}$ and a Hamiltonian diffeomorphism on $\R^{2n}$ that is arbitrarily Hofer-close to $\id$ and has no leafwise fixed points \wrt~$N$.}%
 This may also yield a lower bound on $\big|\Fix(\psi,N)\big|$ that is higher then the sum of the Betti numbers of $N$, for a suitably generic pair $(\psi,N)$. %false if $N$ is closed and regular, see leafwise.tex

Based on the ideas outlined below, one can define a Floer homology for certain regular contact manifolds and use it to show that a given time-1-map of a contact isotopy has translated points. Namely, consider a closed manifold $N$ and a contact form $\al$ on $N$. Assume that $\al$ is regular, i.e., that its symplectic quotient is well-defined. Then $N$ is naturally a smooth principal $S^1$-bundle, which is called a Boothby-Wang fibration. The associated real two-dimensional vector bundle $E$ is equipped with a natural symplectic form, see e.g.~\cite[proof of Lemma 3, p.~200]{GS}. The idea is now to define the Floer homology of $(N,\al)$ to be the Floer homology of the circle bundle in $E$ of some radius $r$. If the symplectic quotient of $N$ is monotone then for a suitable choice of $r$ this circle bundle is a monotone coisotropic submanifold with minimal Maslov number at least 2. Hence its Floer homology is well-defined. Since part of the symplectization of $(N,\al)$ symplectically embeds into $E$, it follows that translated points exist under suitable hypotheses.

Various versions of coisotropic Floer homology may play a role in mirror symmetry, as physicists have realized that the Fukaya category should be enlarged by coisotropic submanifolds, in order to make homological mirror symmetry work, see e.g.~\cite{KO}.

\section*{Floer homology for an adiscal or monotone regular coisotropic submanifold}
To explain the coisotropic Floer homology in the regular case, consider a geometrically bounded symplectic manifold $(M,\om)$, and a closed, %connected
regular coisotropic submanifold $N\sub M$. 

Suppose first also that $N$ is adiscal. We define the \emph{Floer homology of $N$} as follows. Let $\phi=(\phi^t)_{t\in[0,1]}$ be a Hamiltonian isotopy starting at $\id$, such that $(N,\phi^1)$ is nondegenerate. We call a point $x_0\in M$ a \emph{$(N,\phi)$-contractible leafwise fixed point} iff the path $x:[0,1]\to M$, $x(t):=\phi^t(x_0)$, is semistrongly $(N,\om)$-homotopic to a constant path.\footnote{By definition, for every such point $x_0$, the point $x(1)=\phi^1(x_0)$ lies in the isotropic leaf of $x_0$. Hence $x_0$ is a leafwise fixed point of $\phi^1$.} %
We define
\begin{eqnarray}&\nn\Fixc(N,\phi):=\big\{(N,\phi)\textrm{-contractible leafwise fixed points}\big\},&\\
\label{eq:CF N H}&\CF(N,\phi):=\bigoplus_{\Fixc(N,\phi)}\Z_2.&
\end{eqnarray}
\begin{Rmk}By definition this direct sum contains one copy of $\Z_2$ for each point in $\Fixc(N,\phi)$.
\end{Rmk}

We now define a collection of boundary operators on $\CF(N,\phi)$, one for each $(N,\phi^1)$-regular time-dependent almost complex structure. To explain this, observe that, since $N$ is regular, $N_\om$ carries canonical smooth and symplectic structures $\A_{N,\om}$ and $\om_N$. We define
\begin{eqnarray}\nn&\hhat M:=M\x N_\om,\quad\hhat\om:=\om\oplus(-\om_N),&\\
\label{eq:iota N N hhat M}&\iota_N:N\to\hhat M,\,\iota_N(x):=\big(x,\textrm{isotropic leaf through }x\big),\,\hhat N:=\iota_N(N),&\\
\nn&\hhat\phi^t:=\phi^t\x\id_{N_\om}.&
\end{eqnarray} 
The map $\iota_N$ is a $\hhat\om$-Lagrangian embedding of $N$ into $\hhat M$, see \cite[Lemma 3.2]{ZiLeafwise}. %\ref{le:N}
It follows from \cite[Proposition 61, p.~43]{ZiMaslov} that the restriction of the map $\iota_N$, given by
\begin{eqnarray}\label{eq:iota N Fixc}&\iota_N:\Fixc(N,\phi)\to\Fixc(\hhat N,\hhat\phi)=&\\
\nn&\left\{\hhat x\in\hhat N\cap\big(\hhat\phi^1\big)^{-1}(\hhat N)\,\Big|\,t\mapsto\hhat\phi^t(\hhat x)\textrm{ contractible with endpoints in }\hhat N\right\}&
\end{eqnarray}
is well-defined and injective. A straightforward argument shows that it is surjective. %We don't need to assume that $N$ is connected, here.

Let $p\in(2,\infty)$ and $\psi$ be a Hamiltonian diffeomorphism on $M$. We call a $t$-dependent $\hhat\om$-compatible smooth almost complex structure $\hhat J$ on $\hhat M$ \emph{$(N,\psi)$-regular}, iff the vertical differential of its Cauchy-Riemann operator is surjective for all finite-energy $\hhat J$-holomorphic strips with compact image and boundary on $\hhat N$ and $\hhat\psi(\hhat N)$. (For definitions see \cite[Proposition 2.1]{FlLag}.\footnote{The exponent $p$ enters the Banach bundle setup of that proposition.}) We define $\Jreg(N,\psi):=\Jreg\big(M,\om,N,\psi,p\big)$ to be the set of all $(N,\psi)$-regular $\hhat J$.\footnote{It follows from the proof of \cite[Proposition 2.1]{FlLag} that this set is dense in the set of all $t$-dependent $\hhat\om$-compatible almost complex structures, and therefore nonempty. See \cite[Theorem 5]{FlUnreg}.} %end footnote
For every $\hhat J\in\Jreg(N,\phi^1)$ we define
the \emph{Floer boundary operator}
\[\dd_{N,\phi,\hhat J}:\CF(N,\phi)\to\CF(N,\phi)\]
to be the (Lagrangian) Floer boundary operator of $\big(\hhat M,\hhat\om,\hhat N,\hhat\phi,\hhat J\big)$, where $\CF(N,\phi)$ is as in \reff{eq:CF N H}. (See \cite[Definition 3.1]{FlLag}.)

To see that this operator is well-defined, recall that it is defined on the direct sum of $\Z_2$'s indexed by the set occurring on the right hand side of \reff{eq:iota N Fixc}.\footnote{Sometimes this is called the ``$\Z_2$-vector space spanned by this index set''. Here we identify the intersection $\hhat N\cap\big(\hhat\phi^1\big)^{-1}(\hhat N)$ with $\hhat\phi^1(\hhat N)\cap\hhat N$ via the map $\hhat\phi^1$.} %end footnote
Using the bijection \reff{eq:iota N Fixc}, we identify this direct sum with $\CF(N,\phi)$. Hence $\dd_{N,\phi,\hhat J}$ is defined between the right spaces.

We check the conditions of \cite[Definition 3.1]{FlLag}. Since $N$ is closed, the same holds for $N_\om$. Since $(M,\om)$ is geometrically bounded, it follows that $(\hhat M,\hhat\om)$ is geometrically bounded. The Lagrangian $\hhat N$ is diffeomorphic to $N$ via $\iota_N$, and therefore also closed. Since $N$ is adiscal, it follows from \cite[Proposition 61, p.~43]{ZiMaslov} that the same holds for $\hhat N$. Since $\big(N,\phi^1\big)$ is nondegenerate, by \cite[Lemma 3.2(c)]{ZiLeafwise} the intersection of $\hhat N$ and $\hhat\phi^{-1}(\hhat N)$ is transverse. Hence the conditions of \cite[Definition 3.1]{FlLag} are satisfied, and therefore the boundary operator $\dd_{N,\phi,\hhat J}$ is well-defined.\footnote{In \cite{FlLag} Floer assumes that the symplectic manifold is closed. However, the same construction of Floer homology works for geometrically bounded symplectic manifolds. Here we use that we only consider Floer strips with compact image.} %end footnote
By \cite[Lemma 3.2]{FlLag} it squares to 0. Hence we may define the \emph{Floer homology of $\big(N,\phi,\hhat J\big)$} to be the homology
\[\HF\big(N,\phi,\hhat J\big):=H\big(\CF(N,\phi),\dd_{N,\phi,\hhat J}\big).\]
Let $G$ be a $\Z$-grading of the $\Z_2$-vector space $\CF(N,\phi)$, which via the identification \reff{eq:iota N Fixc} is compatible with the Viterbo-Maslov index for strips satisfying Lagrangian boundary conditions.\footnote{By \cite[Proposition 2.4]{FlLag} such a grading exists and each two gradings differ by an additive constant.} %end footnote
Such a $G$ induces a grading on $\HF\big(N,\phi,\hhat J\big)$. For every pair $\hhat J_0,\hhat J_1\in\Jreg(N,\phi^1)$ we denote by
\[\Phi_{\hhat J_0,\hhat J_1}:\HF\big(N,\phi,\hhat J_0\big)\to\HF\big(N,\phi,\hhat J_1\big)\]
the canonical isomorphism provided by the proof of \cite[Proposition 3.1, p.~522]{FlLag}. This isomorphism respects the grading $G$. It does not depend on the choice of $G$. 
\begin{defi}[Floer homology for adiscal coisotropic] We define the \emph{Floer homology of $(N,\phi)$} to be
\[\HF(N,\phi):=\left(\left(\HF\big(N,\phi,\hhat J\big)\right)_{\hhat J\in\Jreg(N,\phi^1)},\left(\Phi_{\hhat J_0,\hhat J_1}\right)_{\hhat J_0,\hhat J_1\in\Jreg(N,\phi^1)}\right).\]
\end{defi}
\begin{Rmks}\begin{itemize}
\item This is a collection of graded $\Z_2$-vector spaces together with grading-preserving isomorphisms.
\item Philosophically, the Floer homology of $(N,\phi)$ is defined to be $\HF\big(N,\phi,\hhat J\big)$, for some choice of $\hhat J$. The collection of isomorphisms $\left(\Phi_{\hhat J_0,\hhat J_1}\right)$ encodes the sense in which this does ``not depend'' on this choice.
\end{itemize}
\end{Rmks}
By the proof of \cite[Theorem 1]{FlLag} $\HF(N,\phi)$ is isomorphic to the singular homology of $\hhat N$ (hence of $N$) with $\Z_2$-coefficients. Since $\HF(N,\phi)$ is generated by the $(N,\phi)$-contractible leafwise fixed points of $\phi^1$, this reproves Theorem \ref{thm:adiscal}.

Suppose now that $N$ is monotone in the sense of \cite{ZiMaslov} and of minimal Maslov number $m(N)\geq2$. %2019 11 11 Here we need to assume this, in order to make the Floer homology well-defined, even though for the application to leafwise fixed points we don't need it, since otherwise the stated bound is 0.
\footnote{We continue to assume that $(M,\om)$ is geometrically bounded, $N$ is closed and regular, and that $(N,\phi^1)$ is nondegenerate.} %footnote
\begin{defi}[Floer homology for monotone coisotropic] We define the \emph{Floer homology of $(N,\phi)$} to be the Floer homology of $(\hhat N,\hhat\phi)$, as defined in \cite{OhFloerLag,OhAdd}.
\end{defi}
Since $N$ is monotone, the same holds for $\hhat N$. The minimal Maslov numbers of $N$ and $\hhat N$ agree, see the proof of \cite[Theorem 3, p.~17]{ZiMaslov}. It follows that the Floer homology of $(N,\phi)$ is well-defined. By \cite[Corollary 2.1]{AlLag} it is isomorphic to singular homology in degrees $i=\dim N-m(N)+2,\ldots,m(N)-2$. It follows that \reff{eq:Fix sum m N} holds. This reproves Theorem \ref{thm:monotone}.

\section*{Local Floer homology for an arbitrary coisotropic submanifold}
Consider now the situation in which $(M,\om)$ is any symplectic manifold, $N$ an arbitrary closed coisotropic submanifold of $M$, and $\phi=(\phi^t)$ a Hamiltonian flow on $M$ whose restriction to $N$ stays ``$C^0$-close'' to the inclusion $N\to M$, such that $(N,\phi^1)$ is nondegenerate. We also fix an $\om$-compatible almost complex structure $J$ on $M$. Heuristically, we define the local Floer homology $\HF(N,\phi,J)$ as follows. Its chain complex is generated by the points $x\in\Fix(N,\phi^1)$, for which there is a ``short'' path from $x$ to $\phi^1(x)$ within the isotropic leaf through $x$. 

To explain the boundary operator $\dd=\dd_{N,\phi,J}$, we denote by $i_N:N\to M$ the inclusion map. We equip the product $M\x N$ with the presymplectic form $\om\oplus(-i_N^*\om)$. By \cite[Lemma 4]{ZiC0} there exists a symplectic submanifold $\wt M$ of $M\x N$ that contains the diagonal
\begin{equation}\label{eq:wt N}\wt N:=\big\{(x,x)\,\big|\,x\in N\big\}\end{equation}
as a Lagrangian submanifold. We shrink $\wt M$, so that it is a Weinstein neighbourhood of $\wt N$. The flow $\phi$ induces a Hamiltonian flow $\wt\phi$ that is defined on an open neighbourhood $\wt N$ of $\wt M$. %between two open neighbourhoods of $\wt N$ in $\wt M$.
The structure $J$ induces an almost complex structure $\wt J$ on $\wt M$ that is $\wt\om$-compatible. 

The boundary operator $\dd$ is now defined to be the boundary operator of the ``local Lagrangian Floer homology'' of $\left(\wt M,\wt\om,\wt N,\wt\phi,\wt J\right)$. This map counts finite energy $\wt J$-holomorphic strips in $\wt M$ that stay ``close'' to $\wt N$, have Viterbo-Maslov index $1$, map the lower and upper boundaries of the strip to $\wt N$ and $(\wt\phi^1)^{-1}(\wt N)$, and connect two intersection points of $\wt N$ and $(\wt\phi^1)^{-1}(\wt N)$. Such points correspond to points $x\in\Fix(N,\phi^1)$, for which there exists a short path from $x$ to $\phi^1(x)$ within a leaf. (See \cite[Lemma 6]{ZiC0}.)

To understand why heuristically, the boundary operator $\dd$ is well-defined and squares to zero, observe that $\wt N$ intersects $(\wt\phi^1)^{-1}(\wt N)$ transversely, since $(N,\phi^1)$ is nondegenerate. (See \cite[Lemma 7]{ZiC0}.) Therefore, for generic $\wt J$%
\footnote{Here one needs to work with a family of almost complex structures depending on the time $t$.}%end footnote
, the moduli space of $\wt J$-strips is a 0-dimensional manifold in a natural sense. (Here we divided by the translation action.) It is compact for the following reasons:
\begin{itemize}
\item Holomorphic strips with boundary on $\wt N$ and $(\wt\phi^1)^{-1}(\wt N)$ stay inside some fixed compact neighbourhood of $\wt N$, provided that $\wt\phi^1$ is close enough to the identity. This follows from the fact that there is a neighbourhood $\wt U$ of $\wt N$ and an exhausting $\wt J$-plurisubharmonic function on $\wt U$.
\item Disks or spheres cannot bubble off. This follows from our assumption that $\wt M$ is a Weinstein neighbourhood of $\wt N$, which implies that $\wt N$ is an exact Lagrangian in $\wt M$.
\item Index-1-strips generically do not break.
\end{itemize}
It follows that heuristically, $\dd$ is well-defined. For similar reasons we have $\dd^2=0$.

Given two choices of symplectic submanifolds $\wt M,\wt M'\sub M\x N$ containing $\wt N$, one obtains a symplectomorphism between open neighbourhoods of $\wt N$ in $\wt M$ and $\wt M'$, by sliding $\wt M$ to $\wt M'$ along the isotropic leaves of $N$. This symplectomorphism intertwines the corresponding $\wt\phi$'s and $\wt J$'s. It follows that the boundary operator does not depend on the choice of $\wt M$, and therefore, heuristically, is well-defined. 

To make the outlined Floer homology rigorous, the words ``close'' and ``short'' used above, need to be made precise. To obtain an object that does not depend on the choice of ``closeness'', the \emph{local Floer homology of $(N,J)$} should really be defined to be the \emph{germ} of the map
\[\phi\mapsto\HF(N,\phi,J)\]
around $\id:M\to M$. 

By showing that $\HF(N,\phi,J)$ is isomorphic to the singular homology of $\wt N$, it should be possible to reproduce the lower bound \reff{eq:Fix sum } of Theorem \ref{thm:C0}.

\begin{Rmk}[local presymplectic Floer homology] A presymplectic form on a manifold is a closed two-form with constant rank. By \cite[Proposition 3.2]{Ma} every presymplectic manifold can be coisotropically embedded into some symplectic manifold. By \cite[4.5.~Th\'eor\`eme on p.~79]{Ma} each two coisotropic embeddings are equivalent. Hence heuristically, we may define the local Floer homology of a presymplectic manifold to be the local Floer homology of any of its coisotropic embeddings.
\end{Rmk}

\begin{Rmk}[relation between the constructions]\label{rmk:regular} Assume that $N$ is regular. Then the constructions of its ``global'' and local Floer homologies are related as follows. Namely, the symplectic submanifold $\wt M\sub M\x N$ can be viewed as a local version of $\hhat M$. More precisely, shrinking $\wt M$ if necessary (so that it still contains $\wt N$), the manifold $\wt M$ symplectically embeds into $\hhat M$ via the map
\[(x,y)\mapsto\big(x,\textrm{isotropic leaf through }y\big).\]
\end{Rmk}

\section*{Acknowledgments}

I would like to thank Will Merry for an interesting discussion and the anonymous referees for valuable suggestions.


\begin{thebibliography}{99}

\bibitem[Al1]{AlFloer} P.~Albers, \emph{A note on local Floer homology}, arXiv:math/0606600. %not published in 2014

\bibitem[Al2]{AlLag} P.~Albers, \emph{A Lagrangian Piunikhin-Salamon-Schwarz morphism and two comparison homomorphisms in Floer homology},  Int.~Math.~Res.~Not.~2008,  {\bf no.~4}, Art.~ID~rnm134. 

%\bibitem[AF1]{AFRab} P.~Albers, U.~Frauenfelder, \emph{Leaf-wise intersections and Rabinowitz Floer homology}, J.~Topol.~Anal.~{\bf 2} (2010), {\bf no.~1}, 77--98. 

%\bibitem[AF2]{AFinf} P.~Albers, U.~Frauenfelder, \emph{Infinitely many leaf-wise intersections on cotangent bundles}, Expo.~Math.~{\bf 30} (2012), {\bf no.~2}, 168--181. 

%\bibitem[AF3]{AFSurvey} P.~Albers, U.~Frauenfelder, \emph{Rabinowitz Floer Homology: a Survey}, Global Differential Geometry, Springer Proceedings in Mathematics Vol.~\textbf{17}, 2012, 437--461. 

%\bibitem[AMc]{AMc} P.~Albers, M.~McLean, \emph{Non-displaceable contact embeddings and infinitely many leaf-wise intersections}, J.~Symplectic Geom.~\textbf{9} (2011), \textbf{no.~3}, 271--284. 

%\bibitem[AMo]{AMo} P.~Albers, A.~Momin, \emph{Cup-length estimates for leaf-wise intersections}, Math.~Proc.~Cambridge Philos.~Soc.~\textbf{149} (2010), \textbf{no.~3}, 539--551. 

%\bibitem[Ar]{Ar} V.~I.~Arnold, \emph{Sur une propri\'et\'e topologique des applications globalement canoniques de la m\'ecanique classique}, C.~R.~Acad.~Sci.~Paris {\bf 261} (1965), 3719--3722. 

%\bibitem[Bae]{Bae} Y.~Bae, \emph{On magnetic leaf-wise intersections}, J.~Topol.~Anal.~\textbf{4} (2012), \textbf{no.~3}, 361--411.

%\bibitem[Ban]{Ban} A.~Banyaga, \emph{On fixed points of symplectic maps}, Invent. Math. {\bf 56} (1980), no. 3, 215--229. 

\bibitem[Ch]{Ch} Yu. V. Chekanov, \emph{Lagrangian intersections, symplectic energy, and areas of holomorphic curves}, Duke Math. J. {\bf 95} (1998), {\bf no. 1}, 213-226. 

\bibitem[CFHW]{CFHW} K.~Cieliebak, A.~Floer, H.~Hofer, K.~Wysocki, \emph{Applications of symplectic homology, II, Stability of the action spectrum}, Math.~Z.~\textbf{223} (1996), \textbf{no.~1}, 27--45. 

%\bibitem[Dr]{Dr} D. Dragnev, \emph{Symplectic rigidity, symplectic fixed points, and global perturbations of Hamiltonian systems}, Comm. Pure Appl. Math. {\bf 61} (2008), {\bf no. 3}, 346--370. 

%\bibitem[EH]{EH} I. Ekeland, H. Hofer, \emph{Two symplectic fixed-point theorems with applications to Hamiltonian dynamics}, J.~Math.~Pures Appl.~(9) {\bf 68} (1989), {\bf no. 4}, 467--489 (1990). 

\bibitem[Fl1]{FlLag} A.~Floer, \emph{Morse theory for Lagrangian intersections}, J. Differential Geom. Vol. 28, no. {\bf 3} (1988), 513-547.

\bibitem[Fl2]{FlUnreg} A.~Floer, \emph{The unregularized gradient flow of the symplectic action}, Comm.~Pure Appl.~Math.{\bf 41} (1988), no.~6, 775--813.

\bibitem[Fl3]{FlFixed} A.~Floer, \emph{Symplectic fixed points and holomorphic spheres}, Comm.~Math.~Phys.~{\bf 120} (1989), no.~4, 575--611.

%\bibitem[Ge]{Geiges} H.~Geiges, \emph{An introduction to contact topology}, Cambridge Studies in Advanced Mathematics, \textbf{109}. Cambridge University Press, Cambridge, 2008. %xvi+440 pp. ISBN: 978-0-521-86585-2 

\bibitem[GS]{GS} H.~Geiges, A.~Stipsicz, \emph{Contact structures on product five-manifolds and fibre sums along circles}, Math.~Ann.~\textbf{348} (2010), no.~\textbf{1}, 195--210. 

%\bibitem[Gi]{Gi} V.~L.~Ginzburg, \emph{Coisotropic intersections}, Duke Math. J. {\bf 140} (2007), {\bf no.1}, 111-163. 

\bibitem[GG1]{GGLocal} V.~L.~Ginzburg, B.~G\"urel, \emph{Local Floer homology and the action gap}, J.~Symplectic Geom.~\textbf{8} (2010), \textbf{no.~3}, 323--357. 

\bibitem[GG2]{GGFrag} V.~L.~Ginzburg, B.~G\"urel, \emph{Fragility and persistence of leafwise intersections}, Math.~Z.~\textbf{280} (2015), \textbf{no.~3-4}, 989--1004.

%\bibitem[Gr]{Gr} M.~Gromov, \emph{Pseudoholomorphic curves in symplectic manifolds}, Invent.~Math.~{\bf 82} (1985), no.~2, 307--347. 

%\bibitem[G\"u]{Gu} B.~G\"urel, \emph{Leafwise coisotropic intersections}, Int.~Math.~Res.~Not.~2010, no.~\textbf{5}, 914--931.

%\bibitem[Ho1]{HoLag} H.~Hofer, \emph{Lagrangian embeddings and critical point theory}, Ann.~Inst.~H.~Poincar\'e Anal.~Non Lin\'eaire 2 (1985), \textbf{no.~6}, 407--462. 

%\bibitem[Ho2]{HoTop} H.~Hofer, \emph{On the topological properties of symplectic maps}, Proc. Roy. Soc. Edinburgh Sect. {\bf A 115} (1990), {\bf no. 1-2}, 25--38.

%\bibitem[Ka1]{KaEx} J.~Kang, \emph{Existence of leafwise intersection points in the unrestricted case}, Israel Journal of Mathematics, vol.~\textbf{190} (2012), iss.~\textbf{1}, 111--134. 

%\bibitem[Ka2]{KaGen} J.~Kang, \emph{Generalized Rabinowitz Floer homology and coisotropic intersections}, Int.~Math.~Res.~Not.~\textbf{no.~10} (2013), 2271--2322. 

\bibitem[KO]{KO} A.~Kapustin and D.~Orlov, \emph{Remarks on $A$-branes, mirror symmetry, and the Fukaya category}, J.~Geom.~Phys.~{\bf 48} (2003), 84--99.

%\bibitem[MMP]{MMP} L.~Macarini, W.~J.~Merry, G.~P.~Paternain, \emph{On the growth rate of leaf-wise intersections}, J.~Symplectic Geom.~\textbf{10} (2012), \textbf{no.~4}, 601--653.

\bibitem[Ma]{Ma} C.-M. Marle, \emph{Sous-vari\'et\'es de rang constant d'une vari\'et\'e symplectique}, Ast\'erisque, {\bf 107-108}, Soc. Math. France, Paris, 1983. 

%\bibitem[Mo]{Mo} J.~Moser, \emph{A fixed point theorem in symplectic geometry}, Acta Math. {\bf 141} (1978), no. 1-2, 17--34.

%keep remark: don't cite, since only small part about local Floer homology 
%\bibitem[Oh1]{OhRelative} Y.-G.~Oh, \emph{Relative Floer and quantum cohomology and the symplectic topology of Lagrangian submanifolds}, Contact and symplectic geometry (Cambridge, 1994), 201--267, Publ.~Newton Inst., \textbf{8}, Cambridge Univ.~Press, Cambridge, 1996. %not available through Mathscinet/ UU

\bibitem[Oh1]{OhFloerLag} Y.-G.~Oh, \emph{Floer cohomology of Lagrangian intersections and pseudo-holomorphic disks. I.} Comm.~Pure Appl.~Math.~\textbf{46} (1993), no.~\textbf{7}, 949--993. 

\bibitem[Oh2]{OhAdd} Y.-G.~Oh, \emph{Addendum to: ``Floer cohomology of Lagrangian intersections and pseudo-holomorphic disks. I.''} Comm.~Pure Appl.~Math.~\textbf{46} (1993), no.~\textbf{7}, 949--993. % MR1223659

\bibitem[Oh3]{OhFloer} Y.-G.~Oh, \emph{Floer cohomology, spectral sequences, and the Maslov class of Lagrangian embeddings}, Internat.~Math.~Res.~Notices 1996, \textbf{no.~7}, 305--346. 

\bibitem[Oh4]{OhLocal} Y.-G.~Oh, \emph{Localization of Floer homology of engulfable topological Hamiltonian loop}, arXiv:1111.5996v4.

\bibitem[Oh5]{OhSympl} Y.-G.~Oh, \emph{Symplectic Topology and Floer Homology}, \url{https://www.math.wisc.edu/~oh/all.pdf}. %keep remark % 2014 08 13 e-mail from Oh: will appear in Cambridge University Press

%\bibitem[Os]{Os} Y.~Ostrover, \emph{A comparison of Hofer s metrics on Hamiltonian diffeomorphisms and Lagrangian submanifolds}, Commun.~Contemp.~Math.~{\bf 5} (2003), {\bf no.~5}, 803--811. 

\bibitem[Po]{Po} M.~Po\'zniak, \emph{Floer homology, Novikov rings and clean intersections}, Northern California Symplectic Geometry Seminar, 119--181, Amer.~Math.~Soc.~Transl.~Ser.~\textbf{2}, \textbf{196}, Amer.~Math.~Soc., Providence, RI, 1999. 

%\bibitem[Sa1]{SaIt} S.~Sandon, \emph{On iterated translated points for contactomorphisms of $\R^{2n=1}$ and $\R^{2n}\x S^1$}, Internat.~J.~Math.~\textbf{23} (2012), \textbf{no.~2}, 1250042, 14 pp.

\bibitem[Sa2]{SaMorse} S.~Sandon, \emph{A Morse estimate for translated points of contactomorphisms of spheres and projective spaces}, Geom.~Dedicata \textbf{165} (2013), 95--110.

%\bibitem[SZ]{SZ} J.~Swoboda and F.~Ziltener, \emph{Coisotropic Displacement and Small Subsets of a Symplectic Manifold}, Math.~Z., Vol.~{\bf 271}, Iss.~{\bf 1} (2012), p.~415--445.

%\bibitem[We1]{WePert} A.~Weinstein, \emph{$C^0$ perturbation theorems for symplectic fixed points and Lagrangian intersections}, South Rhone seminar on geometry, III (Lyon, 1983), 140--144, Travaux en Cours, Hermann, Paris, 1984. 

%\bibitem[We2]{WeCrit} A.~Weinstein, \emph{Critical point theory, symplectic geometry and Hamiltonian systems}, Proceedings of the 1983 Beijing symposium on differential geometry and differential equations, 261--289, Science Press, Beijing, 1986.

\bibitem[Zi1]{ZiLeafwise} F. Ziltener, \emph{Coisotropic Submanifolds, Leafwise Fixed Points, and Presymplectic Embeddings}, J.~Symplectic Geom.~{\bf 8} (2010), no.~{\bf 1}, 1--24.

\bibitem[Zi2]{ZiMaslov} F.~Ziltener, \emph{A Maslov Map for Coisotropic Submanifolds, Leaf-wise Fixed Points and Presymplectic Non-Embeddings}, arXiv:0911460.

\bibitem[Zi3]{ZiC0} F.~Ziltener, \emph{Leafwise fixed points for $C^0$-small Hamiltonian flows}, arXiv:1408.4578, accepted by IMRN.

\end{thebibliography}
\end{document}